\documentclass{article}
\usepackage{amssymb}
\usepackage{amsmath}
\usepackage{graphicx}

\newtheorem{theorem}{Theorem}
\newtheorem{remark}[theorem]{Remark}

\begin{document}

\title{On $p$-adic Twisted Euler $\left( h,q\right) $-$l$-Functions}
\author{S. H. Rim~$^a$, Y. Simsek~$^b$, V. Kurt~$^b$ and T. Kim~$^c$ \\
$^a$Department of Mathematics Education, Kyungpook National University, \\
Daegu 702-701, S. Korea\\
$^b$Department of Mathematics, Akdeniz University, 07058 Antalya,
Turkey\footnote{Second and third authors are supported by Akdeniz
University
Scientific Research Project Unit.}\\
$^c$EECS, Kyungpook National University, Daegu 702-701, S. Korea\\
[0pt] shrim@knu.ac.kr, ysimsek@akdeniz.edu.tr, vkurt@akdeniz.edu.tr,
tkim@knu.ac.kr}
\date{}
\maketitle

\textbf{Abstract : }In the recent paper, Kim-Rim have studied interesting
twisted $q$-Euler numbers and polynomials. In \cite{Kim-Rim2007}, Kim-Rim
suggested the question to find a $q$-analogue of the $p$-adic twisted $%
\left( h,q\right) $-$l$-function which interpolates generalized twisted $%
\left( h,q\right) $-Euler numbers attached to $\chi $. This question is
remained open. The purpose of this paper is to give the answer of the
question.

\bigskip

\section{{\protect\normalsize {Introduction}}}

\hspace{0.2in}Let $p$ be a fixed odd prime number. Throughout this paper, $%
\mathbb{Z}_{p}$, $\mathbb{Q}_{p}$, $\mathbb{C}$ and $\mathbb{C}_{p}$ are
respectively denoted as the ring of $p$-adic rational integers, the field of
$p$-adic rational numbers, the complex numbers field and the completion of
algebraic closure of $\mathbb{Q}_{p}$. Let $v_{p}$ be the normalized
exponential valuation of $\mathbb{C}_{p}$ with $\left| p\right|
_{p}=p^{-v_{p}\left( p\right) }=\frac{1}{p}$.

When one talks of $q$-extension, $q$ is considered in many ways such as an
indeterminate, a complex number $q\in \mathbb{C}$, or $p$-adic number $q\in
\mathbb{C}_{p}$. If $q\in \mathbb{C}$ one normally assumes that $\left|
q\right| <1$. If $q\in \mathbb{C}_{p}$, we normally assume that $\left|
1-q\right| _{p}<p^{-\frac{1}{p-1}}$, so that $q^{x}=$exp$\left( x\text{log}%
q\right) $ for $\left| x\right| _{q}\leqslant 1$.

We use the notations as%
\begin{eqnarray*}
\left[ x\right] _{q} &=&\frac{1-q^{x}}{1-q}=1+q+q^{2}+\cdots +q^{x-1}, \\
\left[ x\right] _{-q} &=&\frac{1-\left( -q\right) ^{x}}{1+q}%
=1-q+q^{2}-q^{3}+\cdots +\left( -1\right) ^{x}q^{x-1}.
\end{eqnarray*}

Let $UD\left( \mathbb{Z}_{p}\right) $ be the set of uniformly differentiable
function on $\mathbb{Z}_{p}$. For $f\in UD\left( \mathbb{Z}_{p}\right) $,
Kim originally defined the $p$-adic invariant $q$-integral on $\mathbb{Z}%
_{p} $ as follows:%
\begin{equation*}
I_{q}\left( f\right) =\int\limits_{\mathbb{Z}_{p}}f\left( x\right) d\mu
_{q}\left( x\right) =\underset{N\rightarrow \infty }{\text{lim}}\frac{1}{%
\left[ p^{N}\right] _{q}}\sum_{x=0}^{p^{N}-1}f\left( x\right) q^{x}\text{,
cf. \cite{KimArxiva}, \cite{Kim2002}, \cite{Simsek2006a},}
\end{equation*}

\noindent where $N$ is natural number. Let%
\begin{equation*}
I_{1}\left( f\right) =\underset{q\rightarrow 1}{\text{lim}}I_{q}\left(
f\right) =\int\limits_{\mathbb{Z}_{p}}f\left( x\right) d\mu _{1}\left(
x\right) =\underset{N\rightarrow \infty }{\text{lim}}\frac{1}{p^{N}}%
\sum_{x=0}^{p^{N}-1}f\left( x\right) ,
\end{equation*}

\noindent where $N$ is a natural number (see \cite{Cenkci-Can2006}, \cite%
{Kim2006a}, \cite{Kim2006b}, \cite{Kim2007}, \cite{Kim-Jang-Rim-Pak2003}).

Let $d$ be a fixed integer. For any positive integer $N$, we set%
\begin{eqnarray*}
\mathbb{X} &=&\mathbb{X}_{d}=\underset{N}{\underleftarrow{\text{lim}}}\left(
\mathbb{Z}/dp^{N}\mathbb{Z}\right) \text{,} \\
\mathbb{X}^{\ast } &=&\bigcup\limits_{\underset{\left( a,p\right) =1}{0<a<dp}%
}\left( a+dp\mathbb{Z}_{p}\right) , \\
a+dp^{N} &=&\left\{ x\in \mathbb{X}:x\equiv a\left( \text{mod}dp^{N}\right)
\right\} \text{,}
\end{eqnarray*}

where $a\in \mathbb{Z}$ lies in $0\leqslant a<dp^{N}$.

Let us define $I_{-q}\left( f\right) $ as%
\begin{equation*}
I_{-q}\left( f\right) =\underset{q\rightarrow -q}{\text{lim}}I_{q}\left(
f\right) =\int\limits_{\mathbb{Z}_{p}}f\left( x\right) d_{-q}\left( x\right)
..
\end{equation*}

\noindent This integral, $I_{-q}\left( f\right) $, can be considered as the $%
q$-deformed $p$-adic invariant integral on $\mathbb{Z}_{p}$ in the sense of
fermionic, cf. \cite{Kim2006a}, \cite{Kim2006b}, \cite{Kim2007}, \cite%
{KimArxiva}, \cite{Kim-Rim2007}.

In \cite{KimArxiva}, multiple $q$-Euler polynomials of higher order were
defined by%
\begin{equation*}
E_{n,q}^{\left( h,k\right) }\left( x\right) =\underset{k\text{-times}}{%
\int\limits_{\mathbb{Z}_{p}}\cdots \int\limits_{\mathbb{Z}_{p}}}\left[
x+x_{1}+\cdots +x_{k}\right] _{q}^{n}q^{\sum_{i=1}^{k}x_{i}\left( h-i\right)
}d\mu _{-q}\left( x\right) \cdots d\mu _{-q}\left( x\right) ,
\end{equation*}

\noindent where $h\in \mathbb{Z}$, $k\in \mathbb{N}$. The $q$-Euler
polynomials of higher order at $x=0$ are called $q$-Euler numbers of higher
order.

In \cite{Carlitz1954}, Carlitz originally constructed $q$-Bernoulli numbers
and polynomials. These numbers and polynomials are studied by many authors
(see \cite{Kim2002}, \cite{Kim-Rim2007}, \cite{Simsek2006a}, \cite%
{Simsek2006b}). In particular, twisted $\left( h,q\right) $-Bernoulli
numbers and polynomials were also studied by several authors (see \cite%
{Kim-Jang-Rim-Pak2003}, \cite{Simsek2006a}, \cite{Simsek2006b}).

In \cite{Kim-Rim2007}, Kim-Rim introduced an interesting twisted $q$-Euler
numbers and polynomials associated with basic twisted $q$-$l$-functions and
suggested the following question:

``Find a $q$-analogue of the $p$-adic twisted $l$-function which
interpolates generalized twisted $q$-Euler numbers attached to $\chi $, $%
E_{n,w,\chi ,q}$.''

In this paper, we give some interesting properties related to twisted $%
\left( h,q\right) $-Euler numbers and polynomials. The final purpose of this
paper is to construct $p$-adic twisted Euler $\left( h,q\right) $-$l$%
-function which is a part of answer for the question in \cite{Kim-Rim2007}.

\section{{\protect\normalsize {$p$-adic Invariant Integral on $\mathbb{Z}_{p}
$ Associated with Twisted $\left( h,q\right) $-Euler Numbers and Polynomials}%
}}

\hspace{0.2in}Let $h\in \mathbb{Z}$ and $q\in \mathbb{C}_{p}$ with $\left|
1-q\right| _{p}<p^{-\frac{1}{p-1}}$. From the invariant integral on $\mathbb{%
Z}_{p}$ in the sense of fermionic, we define%
\begin{equation*}
I_{-1}\left( f\right) =\underset{q\rightarrow -1}{\text{lim}}I_{q}\left(
f\right) =\int\limits_{\mathbb{Z}_{p}}f\left( x\right) d\mu _{-1}\left(
x\right) ,
\end{equation*}

\noindent where $f\in UD\left( \mathbb{Z}_{p}\right) $, cf. \cite{Kim2006a}.
Note that $I_{-1}\left( f_{1}\right) +I_{-1}\left( f\right) =2f\left(
0\right) $, where $f_{1}\left( x\right) =f\left( x+1\right) $. Let $%
C_{p^{n}}=\left\{ \xi :\xi ^{p^{n}}=1\right\} $ be the cyclic group of order
$p^{n}$ and let $T_{p}=\underset{n\rightarrow \infty }{\text{lim}}%
C_{p^{n}}=C_{p^{\infty }}$. Then $T_{p}$ is $p$-adic locally constant space.
For $\xi \in T_{p}$, we denote by $\phi _{\xi }:\mathbb{Z}_{p}\rightarrow
\mathbb{C}_{p}$ defined by $\phi _{\xi }\left( x\right) =\xi ^{x}$ be the
locally constant function. If we take $f\left( x\right) =\phi _{\xi }\left(
x\right) e^{tx}$, then we have%
\begin{equation}
\int\limits_{\mathbb{Z}_{p}}e^{tx}\phi _{\xi }\left( x\right) d\mu
_{-1}\left( x\right) =\frac{2}{\xi e^{t}+1}\text{, (see \cite{Kim2006a}).}
\label{1}
\end{equation}

\noindent In complex case the twisted Euler numbers were defined by Kim-Rim %
\cite{Kim-Rim2007}%
\begin{equation*}
\frac{2}{\xi e^{t}+1}=\sum_{n=0}^{\infty }E_{n,\xi }\frac{t^{n}}{n!},
\end{equation*}

\noindent where $\left| \text{log}\xi +t\right| <\pi $. Thus we have%
\begin{equation*}
\int\limits_{\mathbb{Z}_{p}}x^{n}\phi _{\xi }\left( x\right) d\mu
_{-1}\left( x\right) =E_{n,\xi }\text{, }\left( n\geqslant 0\right) .
\end{equation*}

\noindent By using iterative method of $p$-adic invariant integral on $%
\mathbb{Z}_{p}$ in the sense of fermionic, we see that%
\begin{equation}
I_{-1}\left( f_{n}\right) =\left( -1\right) ^{n}I_{-1}\left( f\right)
+2\sum_{l=0}^{n-1}\left( -1\right) ^{n-1-l}f\left( l\right) ,  \label{2}
\end{equation}

\noindent where $f_{n}\left( x\right) =f\left( x+n\right) $. If $n$ is odd
positive integer, then we have%
\begin{equation}
I_{-1}\left( f_{n}\right) +I_{-1}\left( f\right) =2\sum_{l=0}^{n-1}\left(
-1\right) ^{l}f\left( l\right) .  \label{3}
\end{equation}

Let $\chi $ be the Dirichlet's character with conductor $d$ ($=$odd)$\in
\mathbb{N}$, and let $f\left( x\right) =\chi \left( x\right) \phi _{\xi
}\left( x\right) e^{tx}\in UD\left( \mathbb{Z}_{p}\right) $. From (\ref{3}),
we can derive the following:%
\begin{equation}
\int\limits_{\mathbb{X}}e^{tx}\phi _{\xi }\left( x\right) \chi \left(
x\right) d\mu _{-1}\left( x\right) =2\frac{\sum_{a=0}^{d-1}\chi \left(
a\right) \phi _{\xi }\left( x\right) e^{at}}{\xi ^{d}e^{dt}+1}.  \label{4}
\end{equation}

\noindent Now we define twisted generalized Euler numbers attached to $\chi $
as follows:%
\begin{equation*}
\frac{2\sum_{a=0}^{d-1}\chi \left( a\right) \phi _{\xi }\left( x\right)
e^{at}}{\xi ^{d}e^{dt}+1}=\sum_{n=0}^{\infty }E_{n,\chi ,\xi }\frac{t^{n}}{n!%
}.
\end{equation*}

\noindent From (\ref{1}) and (\ref{4}), we note that%
\begin{eqnarray}
\int\limits_{\mathbb{Z}_{p}}x^{n}\phi _{\xi }\left( x\right) d\mu
_{-1}\left( x\right) &=&E_{n,\xi },  \notag \\
\int\limits_{\mathbb{X}}x^{n}\phi _{\xi }\left( x\right) \chi \left(
x\right) d\mu _{-1}\left( x\right) &=&E_{n,\chi ,\xi }  \label{5}
\end{eqnarray}

In \cite{KimArxiva}, $\left( h,q\right) $-Euler numbers were defined by%
\begin{equation*}
E_{n,q}^{\left( h,1\right) }\left( x\right) =\int\limits_{\mathbb{Z}%
_{p}}q^{\left( h-1\right) y}\left[ x+y\right] _{q}^{n}d\mu _{-q}\left(
y\right) ,
\end{equation*}

\noindent where $h\in \mathbb{Z}$. $E_{n,q}^{\left( h,1\right) }\left(
0\right) =E_{n,q}^{\left( h,1\right) }$ will be called $\left( h,q\right) $%
-Euler number.\ In the special case $h=1$, we note that $\underset{%
h\rightarrow 1}{\text{lim}}E_{n,q}^{\left( h,1\right) }=E_{n,q}^{\left(
1,1\right) }$ becomes $q$-Euler numbers, which were originally defined by
Kim \cite{KimArxiva}. That is $\underset{h\rightarrow 1}{\text{lim}}%
E_{n,q}^{\left( h,1\right) }=E_{n,q}^{\left( 1,1\right) }=E_{n,q}$.

In the viewpoint of (\ref{5}), we consider twisted $\left( h,q\right) $%
-Euler numbers using $p$-adic invariant $q$-integral on $\mathbb{Z}_{p}$ in
the sense of fermionic as follows:%
\begin{equation}
E_{n,\xi ,q}^{\left( h,1\right) }\left( x\right) =\int\limits_{\mathbb{Z}%
_{p}}q^{\left( h-1\right) y}\phi _{\xi }\left( y\right) \left[ x+y\right]
_{q}^{n}d\mu _{-q}\left( y\right) ,  \label{6}
\end{equation}

\noindent which are called twisted $\left( h,q\right) $-Euler polynomials.
In the special case $x=0$, we use notation $E_{n,\xi ,q}^{\left( h,1\right)
}\left( 0\right) =E_{n,\xi ,q}^{\left( h,1\right) }$ which are called $%
\left( h,q\right) $-twisted Euler numbers. Note that%
\begin{equation}
E_{n,\xi ,q}^{\left( h,1\right) }\left( x\right) =\frac{\left[ 2\right] _{q}%
}{\left( 1-q\right) ^{n}}\sum_{j=0}^{n}\binom{n}{j}\left( -1\right)
^{j}q^{xj}\frac{1}{1+\xi q^{h+j}},  \label{7}
\end{equation}

\noindent where%
\begin{equation*}
\binom{n}{j}=\frac{n\left( n-1\right) \cdots \left( n-j+1\right) }{j!}.
\end{equation*}

\noindent From (\ref{7}), we note that%
\begin{equation}
E_{n,\xi ,q}^{\left( h,1\right) }\left( x\right) =\left[ 2\right]
_{q}\sum_{k=0}^{\infty }\left( -1\right) ^{k}\xi ^{k}q^{hk}\left[ x+k\right]
_{q}^{n},  \label{8}
\end{equation}

\noindent where $h\in \mathbb{Z}$, $n\in \mathbb{N}$. Equation (\ref{8}) is
equivalent to%
\begin{equation*}
E_{n,\xi ,q}^{\left( h,1\right) }\left( x\right) =\frac{\left[ 2\right] _{q}%
}{\left[ 2\right] _{q^{d}}}\left[ d\right] _{q}^{n}\sum_{a=0}^{d-1}\left(
-1\right) ^{a}\xi ^{a}q^{ha}E_{n,\xi ^{d},q^{d}}^{\left( h,1\right) }\left(
\frac{x+a}{d}\right)
\end{equation*}

\noindent (distribution for $E_{n,\xi ,q}^{\left( h,1\right) }\left(
x\right) $), where $n$, $d\left( =\text{odd}\right) \in \mathbb{N}$.

Let $\chi $ be the Dirichlet character with conductor $f\left( =\text{odd}%
\right) \in \mathbb{N}$. Then we define the generalized twisted $\left(
h,q\right) $-Euler numbers attached to $\chi $ as follows: For $n\geqslant 0$%
,%
\begin{equation}
E_{n,\xi ,\chi ,q}^{\left( h,1\right) }=\int\limits_{\mathbb{X}}\chi \left(
x\right) q^{\left( h-1\right) x}\xi ^{x}\left[ x\right] _{q}^{n}d\mu
_{-q}\left( x\right) ,  \label{9}
\end{equation}

\noindent where $h\in \mathbb{Z}$. Note that $E_{n,1,\chi ,q}^{\left(
1,1\right) }=E_{n,\chi ,q}$, see \cite{Kim2006a}. From (\ref{9}), we also
derive%
\begin{eqnarray}
E_{n,\xi ,\chi ,q}^{\left( h,1\right) } &=&\left[ f\right] _{q}^{n}\frac{%
\left[ 2\right] _{q}}{\left[ 2\right] _{q^{f}}}\sum_{a=0}^{f-1}\chi \left(
a\right) \left( -1\right) ^{a}\xi ^{a}q^{ha}E_{n,\xi ^{f},q^{f}}^{\left(
h,1\right) }\left( \frac{a}{f}\right)   \notag \\
&=&\left[ 2\right] _{q}\sum_{k=1}^{\infty }\chi \left( k\right) \left(
-1\right) ^{k}\xi ^{k}q^{hk}\left[ k\right] _{q}^{n},  \label{10}
\end{eqnarray}

\noindent where $n$, $d\left( =\text{odd}\right) \in \mathbb{N}$.

\section{{\protect\normalsize {Twisted $\left( h,q\right) $-Euler Zeta
Function in $\mathbb{C}$}}}

\hspace{0.2in}For $q\in \mathbb{C}$ with $\left| q\right| <1$, $s\in \mathbb{%
C}$, we define%
\begin{equation*}
\zeta _{E,q,\xi }^{\left( h,1\right) }\left( s\right) =\left[ 2\right]
_{q}\sum_{k=1}^{\infty }\frac{\left( -1\right) ^{k}\xi ^{k}q^{hk}}{\left[ k%
\right] _{q}^{s}}.
\end{equation*}

\noindent Then we see that $\zeta _{E,q,\xi }^{\left( h,1\right) }\left(
s\right) $ is analytic continuation in whole complex plane. We easily see
that%
\begin{equation*}
\zeta _{E,q,1}^{\left( h,1\right) }\left( s\right) =\zeta _{E,q}^{\left(
h,1\right) }\left( s\right) \text{, see \cite{KimArxiva}.}
\end{equation*}

We now also consider Hurwitz's type twisted $\left( h,q\right) $-Euler zeta
function as follows: For $s\in \mathbb{C}$, define%
\begin{equation}
\zeta _{E,q,\xi }^{\left( h,1\right) }\left( s,x\right) =\left[ 2\right]
_{q}\sum_{k=0}^{\infty }\frac{\left( -1\right) ^{a}\xi ^{k}q^{hk}}{\left[ x+k%
\right] _{q}^{s}},  \label{11}
\end{equation}

\noindent for $s\in \mathbb{C}$, $h\in \mathbb{Z}$. By (\ref{8}) and (\ref%
{11}), we see that%
\begin{equation*}
\zeta _{E,q,\xi }^{\left( h,1\right) }\left( -n,x\right) =E_{n,\xi
,q}^{\left( h,1\right) }\left( x\right) ,
\end{equation*}

\noindent where $n\in \mathbb{N}$, $h\in \mathbb{Z}$.

Let $\chi $ be the Dirichlet's character with conductor $f\left( =\text{odd}%
\right) \in \mathbb{N}$. Then we define twisted $\left( h,q\right) $-$l$%
-function which interpolates twisted generalized $\left( h,q\right) $-Euler
numbers attached to $\chi $ as follows: For $s\in \mathbb{C}$, $h\in \mathbb{%
Z}$, we define%
\begin{equation}
l_{q,\xi }^{\left( h,1\right) }\left( s,\chi \right) =\left[ 2\right]
_{q}\sum_{k=1}^{\infty }\frac{\chi \left( k\right) \left( -1\right)
^{k}q^{hk}\xi ^{k}}{\left[ k\right] _{q}^{s}}.  \label{12}
\end{equation}

\noindent For any positive integer $n$, we have%
\begin{equation}
l_{q,\xi }^{\left( h,1\right) }\left( -n,\chi \right) =E_{n,\xi ,\chi
,q}^{\left( h,1\right) }\text{, }n\in \mathbb{N}\text{.}  \label{12-1}
\end{equation}

\noindent From (\ref{12-1}), we derive%
\begin{eqnarray*}
l_{q,\xi }^{\left( h,1\right) }\left( s,\chi \right) &=&\left[ 2\right]
_{q}\sum_{k=1}^{\infty }\frac{\chi \left( k\right) \left( -1\right)
^{k}q^{hk}\xi ^{k}}{\left[ k\right] _{q}^{s}} \\
&=&\left[ f\right] _{q}^{-s}\frac{\left[ 2\right] _{q}}{\left[ 2\right]
_{q^{f}}}\sum_{a=1}^{f}\chi \left( a\right) \left( -1\right) ^{a}\xi
^{a}q^{ha}\zeta _{E,\xi ^{f},q^{f}}^{\left( h,1\right) }\left( s,\frac{a}{f}%
\right) .
\end{eqnarray*}

Let $s$ be a complex variable and let $a$ and $F\left( =\text{odd}\right) $
be integer with $0<a<F$. We consider the following twisted $\left(
h,q\right) $-harmonic sums (or partial $\left( h,q\right) $-zeta function):%
\begin{eqnarray*}
H_{E,q,\xi }^{\left( h,1\right) }\left( s,a|F\right)  &=&\sum_{\underset{m>0}%
{m\equiv 0\left( \text{mod}F\right) }}\frac{\left( -1\right) ^{m}q^{hm}\xi
^{m}}{\left[ m\right] _{q}^{s}} \\
&=&\sum_{n=0}^{\infty }\frac{\left( -1\right) ^{a+nF}q^{h\left( a+nF\right)
}\xi ^{a+nF}}{\left[ a+nF\right] _{q}^{s}} \\
&=&\left( -1\right) ^{a}q^{ha}\xi ^{a}\sum_{n=0}^{\infty }\frac{\left(
-1\right) ^{n}\left( q^{F}\right) ^{hn}\left( \xi ^{F}\right) ^{n}}{\left[ F%
\right] _{q}^{s}\left[ \frac{a}{F}+n\right] _{q^{F}}^{s}} \\
&=&\left[ F\right] _{q}^{-s}\frac{\left( -1\right) ^{a}q^{ha}\xi ^{a}}{\left[
2\right] _{q^{F}}}\zeta _{E,\xi ^{F},q^{F}}^{\left( h,1\right) }\left( s,%
\frac{a}{F}\right) .
\end{eqnarray*}

\noindent Thus we have%
\begin{equation}
H_{E,q,\xi }^{\left( h,1\right) }\left( s,a|F\right) =\left[ F\right]
_{q}^{-s}\frac{\left( -1\right) ^{a}q^{ha}\xi ^{a}}{\left[ 2\right] _{q^{F}}}%
\zeta _{E,\xi ^{F},q^{F}}^{\left( h,1\right) }\left( s,\frac{a}{F}\right) .
\label{13}
\end{equation}

\noindent By (\ref{12}), (\ref{12-1}) and (\ref{13}), we see that%
\begin{equation}
l_{q,\xi }^{\left( h,1\right) }\left( s,\chi \right) =\left[ 2\right]
_{q}\sum_{a=1}^{F}\chi \left( a\right) H_{E,q,\xi }^{\left( h,1\right)
}\left( s,a|F\right) .  \label{14}
\end{equation}

\noindent From (\ref{13}), we note that%
\begin{equation}
H_{E,q,\xi }^{\left( h,1\right) }\left( -n,a|F\right) =\frac{\left[ F\right]
_{q}^{n}}{\left[ 2\right] _{q^{F}}}\left( -1\right) ^{a}q^{ha}\xi
^{a}E_{n,\xi ^{F},q^{F}}^{\left( h,1\right) }\left( \frac{a}{F}\right) ,
\label{15}
\end{equation}

\noindent where $n$ is a positive integer. By (\ref{10}), (\ref{14}) and (%
\ref{15}), we see that%
\begin{equation*}
l_{q,\xi }^{\left( h,1\right) }\left( -n,\chi \right) =E_{n,\xi ,\chi
,q}^{\left( h,1\right) }.
\end{equation*}

\noindent The Euler $\left( h,q\right) $-twisted harmonic sum $H_{E,q,\xi
}^{\left( h,1\right) }\left( s,a|F\right) $ will be called partial twisted $%
\left( h,q\right) $-zeta function which interpolates twisted $\left(
h,q\right) $-Euler polynomials at negative integers. The values $l_{q,\xi
}^{\left( h,1\right) }\left( s,\chi \right) $ are algebraic, hence regarded
as lying in the extension of $\mathbb{Q}_{p}.$

\section{{\protect\normalsize {$p$-adic Twisted Euler $\left( h,q\right) $-$l
$-Function}}}

\hspace{0.2in}Let $\omega \left( x\right) $ be the Teichm\"{u}ller character
and let $\left\langle x\right\rangle _{q}=\left\langle x\right\rangle =\frac{%
\left[ x\right] _{q}}{\omega \left( x\right) }$. When $F\left( =\text{odd}%
\right) $ is multiple of $p$, and $\left( a,p\right) =1$, we define $p$-adic
partial $\left( h,q\right) $-zeta function as follows: For $h\in \mathbb{Z}$%
, $q\in \mathbb{C}_{p}$ with $\left| 1-q\right| _{p}<p^{-\frac{1}{p-1}}$, we
define%
\begin{equation*}
H_{E,p,q,\xi }^{\left( h,1\right) }\left( s,a|F\right) =\frac{\left(
-1\right) ^{a}\xi ^{a}q^{ha}}{\left[ 2\right] _{q^{F}}}\left\langle
a\right\rangle ^{-s}\sum_{j=0}^{\infty }\binom{-s}{j}\left( \frac{\left[ F%
\right] _{q}}{\left[ a\right] _{q}}\right) ^{j}q^{aj}E_{j,\xi
^{F},q^{F}}^{\left( h,1\right) },
\end{equation*}

\noindent where $s\in \mathbb{Z}_{p}$. Thus we note that%
\begin{eqnarray*}
H_{E,p,q,\xi }^{\left( h,1\right) }\left( -n,a|F\right) &=&\frac{\left(
-1\right) ^{a}\xi ^{a}q^{ha}}{\left[ 2\right] _{q^{F}}}\left\langle
a\right\rangle ^{n}\sum_{j=0}^{n}\binom{n}{j}\left( \frac{\left[ F\right]
_{q}}{\left[ a\right] _{q}}\right) ^{j}q^{aj}E_{j,\xi ^{F},q^{F}}^{\left(
h,1\right) } \\
&=&\frac{\left( -1\right) ^{a}\xi ^{a}q^{ha}}{\left[ 2\right] _{q^{F}}}\left[
F\right] _{q}^{n}\omega ^{-n}\left( a\right) \sum_{j=0}^{n}\binom{n}{j}%
\left( \frac{\left[ a\right] _{q}}{\left[ F\right] _{q}}\right) ^{n-j}\left(
q^{F}\right) ^{\frac{a}{F}j}E_{j,\xi ^{F},q^{F}}^{\left( h,1\right) } \\
&=&\left[ F\right] _{q}^{n}\frac{\left( -1\right) ^{a}\xi ^{a}q^{ha}}{\left[
2\right] _{q^{F}}}\omega ^{-n}\left( a\right) E_{n,\xi ^{F},q^{F}}^{\left(
h,1\right) }\left( \frac{a}{F}\right) \\
&=&\omega ^{-n}\left( a\right) H_{E,q,\xi }^{\left( h,1\right) }\left(
-n,a|F\right) .
\end{eqnarray*}

\noindent Therefore, we obtain the following formula:%
\begin{equation}
H_{E,p,q,\xi }^{\left( h,1\right) }\left( -n,a|F\right) =\omega ^{-n}\left(
a\right) H_{E,q,\xi }^{\left( h,1\right) }\left( -n,a|F\right) .  \label{16}
\end{equation}

Now we consider $p$-adic interpolating function for twisted generalized $%
\left( h,q\right) $-Euler numbers attached to $\chi $ as follows:%
\begin{equation}
l_{p,q,\xi }^{\left( h,1\right) }\left( s,\chi \right) =\left[ 2\right]
_{q}\sum_{\underset{\left( a,p\right) =1}{a=1}}^{F}\chi \left( a\right)
H_{E,p,q,\xi }^{\left( h,1\right) }\left( s,a|F\right) ,  \label{17}
\end{equation}

\noindent for $s\in \mathbb{Z}_{p}$. Let $n$ be a natural number. Then we
have%
\begin{eqnarray*}
l_{p,q,\xi }^{\left( h,1\right) }\left( -n,\chi \right)  &=&\left[ 2\right]
_{q}\sum_{\underset{\left( a,p\right) =1}{a=1}}^{F}\chi \left( a\right)
H_{E,p,q,\xi }^{\left( h,1\right) }\left( -n,a|F\right)  \\
&=&\left[ 2\right] _{q}\sum_{\underset{\left( a,p\right) =1}{a=1}}^{F}\chi
\left( a\right) \left[ F\right] _{q}^{n}\frac{\left( -1\right) ^{a}\xi
^{a}q^{ha}}{\left[ 2\right] _{q^{F}}}\omega ^{-n}\left( a\right) E_{n,\xi
^{F},q^{F}}^{\left( h,1\right) }\left( \frac{a}{F}\right)  \\
&=&\left[ F\right] _{q}^{n}\frac{\left[ 2\right] _{q}}{\left[ 2\right]
_{q^{F}}}\sum_{\underset{\left( a,p\right) =1}{a=1}}^{F}\left( -1\right)
^{a}\chi \omega ^{-n}\left( a\right) \xi ^{a}q^{ha}E_{n,\xi
^{F},q^{F}}^{\left( h,1\right) }\left( \frac{a}{F}\right)  \\
&=&E_{n,\xi ,\chi \omega ^{-n},q}^{\left( h,1\right) }-\chi \omega
^{-n}\left( p\right) \left[ p\right] _{q}^{n}\frac{\left[ 2\right] _{q}}{%
\left[ 2\right] _{q^{p}}}E_{n,\xi ^{p},\chi \omega ^{-n},q^{p}}^{\left(
h,1\right) }.
\end{eqnarray*}

\noindent Therefore we obtain the following theorem:

\begin{theorem}
\label{thm1}For $s\in \mathbb{Z}_{p}$, we define $p$-adic twisted Euler $%
\left( h,q\right) $-$l$-function as follows:%
\begin{eqnarray*}
l_{p,q,\xi }^{\left( h,1\right) }\left( s,\chi \right) &=&\left[ 2\right]
_{q}\sum_{\underset{\left( a,p\right) =1}{a=1}}^{F}\chi \left( a\right)
H_{E,p,q,\xi }^{\left( h,1\right) }\left( s,a|F\right) \\
&=&\left[ 2\right] _{q}\sum_{\underset{\left( k,p\right) =1}{k=1}}^{\infty }%
\frac{\left( -1\right) ^{k}\chi \left( k\right) q^{hk}\xi ^{k}}{\left[ k%
\right] _{q}^{s}}.
\end{eqnarray*}

\noindent Then we have%
\begin{equation*}
l_{p,q,\xi }^{\left( h,1\right) }\left( -n,\chi \right) =E_{n,\xi ,\chi
\omega ^{-n},q}^{\left( h,1\right) }-\chi \omega ^{-n}\left( p\right) \left[
p\right] _{q}^{n}\frac{\left[ 2\right] _{q}}{\left[ 2\right] _{q^{p}}}%
E_{n,\xi ^{p},\chi \omega ^{-n},q^{p}}^{\left( h,1\right) }.
\end{equation*}
\end{theorem}

\begin{remark}
\label{rem1}From the above theorem, we note that%
\begin{equation*}
l_{p,q,\xi }^{\left( h,1\right) }\left( s,\chi \right) =\int\limits_{\mathbb{%
X}^{\ast }}\chi \left( x\right) \left\langle x\right\rangle ^{-s}q^{\left(
h-1\right) x}\xi ^{x}d\mu _{-q}\left( x\right) .
\end{equation*}
\end{remark}

\end{document}